\documentclass[10pt]{article}
\title{Chromatic Number and Hamiltonicity of Graphs}
\author {Rao Li \\                
         Dept. of mathematical sciences \\
         University of South Carolina Aiken \\
	     Aiken, SC 29801 \\
         {\it Email: raol@usca.edu }
         }
\date{submitted June 11, 2018; accepted Jan. 14, 2019}
\begin{document}
\maketitle
\begin{abstract}
Let $G$ be a $k$ - connected ($k \geq 2$) graph of order $n$.
If $\chi(G) \geq n - k$, then $G$ is Hamiltonian or 
$K_k \vee (K_k^c \cup K_{n - 2k})$ with $n \geq 2 k + 1$,
where $\chi(G)$ is the chromatic number of the graph $G$.  
  \end{abstract} 
$$2010 \,\, Mathematics \,\, Subject \,\, Classification:  05C45, \, 05C15$$
$$Keywords: Hamiltonicity, \,\,  chromatic \,\, number$$ \\ 


\noindent {\bf 1.  Introduction} \\

We consider only finite undirected graphs without loops or multiple edges.
Notation and terminology not defined here follow those in \cite{Bondy}.
Let $G$ be a graph. We use $G^c$ to denote the complement of $G$.
 We also use $\chi(G)$, $\omega(G)$, and $\alpha(G)$ to denote the chromatic number,
 the clique number, and the independent (or stability) number of $G$, respectively. 
We use $G \vee H$ to denote the the join of two disjoint graphs $G$ and $H$. 
If $C$ is a cycle of $G$, we use $\overrightarrow{C}$ to denote the cycle $C$ with a given direction. 
For two vertices $x$, $y$ in $C$, we use $\overrightarrow{C}[x, y]$ to denote the consecutive 
vertices on $C$ from $x$ to $y$ in the direction
specified by $\overrightarrow{C}$. 
The same vertices, in reverse order, are given by $\overleftarrow{C}[y, x]$. 
We use $x^+$ and $x^-$ to denote respectively the successor and predecessor of a vertex $x$ on $C$ along the direction of $C$. 
We also use $x^{++}$ to denote $(x^+)^+$. A cycle $C$  in a graph $G$ is called a Hamiltonian cycle of $G$ if $C$ contains all the vertices of $G$.  
A graph $G$ is called Hamiltonian if $G$ has a Hamiltonian cycle. \\

In this note, we will present a sufficient condition based on the chromatic number 
for the Hamiltonicity of graphs. The main result is as follows. \\
 
 \noindent {\bf Theorem 1.} Let $G$ be a $k$ - connected ($k \geq 2$) graph of order $n$.
If $\chi(G) \geq n - k$, then $G$ is Hamiltonian or
$K_k \vee (K_k^c \cup K_{n - 2k})$ with $n \geq 2 k + 1$.   \\

\noindent {\bf 2.  The Lemmas} \\

We will use the following results as our lemmas. The first one is an inequality established by Nordhaus and Gaddum in \cite{NG}.\\
 
 \noindent {\bf Lemma 1.} Let $G$ be a graph of order $n$. Then
$\chi(G) + \chi(G^c) \leq n + 1$. \\
 
 The second one is the main result in \cite{AM}. \\ 
 
   \noindent {\bf Lemma 2.} Let $G$ be a $k$ - connected ($k \geq 2$) graph with independent number $\alpha = k + 1$.
Let $C$ be the longest cycle in $G$. Then $G[V(G) - V(C)]$ is complete.   \\
 
 \noindent {\bf 3.  Proofs} \\

\noindent{\bf Proof of Theorem 1.} Let $G$ be a $k$ - connected ($k \geq 2$) graph satisfying the conditions in Theorem $1$. 
Assume that $G$ is not Hamiltonian. Then $n \geq 2 k + 1$ (otherwise $\delta \geq k \geq \frac{n}{2}$ and $G$ is Hamiltonian).
Since $k \geq 2$, $G$ contains a cycle. 
Choose a longest cycle $C$ in $G$ and give a direction on $C$. Since $G$ is not Hamiltonian, 
there exists a vertex $x_0 \in V(G) - V(C)$. By Menger's theorem, 
we can find $s$ ($s \geq k$) pairwise disjoint (except for $x_0$) paths $P_1$, $P_2$, ..., $P_s$ 
between $x_0$ and $V(C)$. Let $u_i$ be the end vertex of $P_i$ on $C$, where $1 \leq i \leq s$.  
We assume that the appearance of $u_1$, $u_2$, ..., $u_s$ agrees with the given direction on $C$.
We use $u_i^+$ to denote the successor of $u_i$ along the direction of $C$, where $1 \leq i \leq s$. 
Then a standard proof in Hamiltonian graph theory yields that $T := \{x_0, u_1^+, u_2^+, ..., u_s^+ \}$ 
is independent (otherwise $G$ would have cycles which are longer than $C$). Since $s \geq k$,
we have an independent set $S := \{x_0, u_1^+, u_2^+, ..., u_k^+ \}$ of size $k + 1$ in $G$ and a clique $S$ of size $k + 1$ in $G^c$. 
From Lemma $1$, we have that
$$n + 1 = n - k + k + 1 \leq \chi(G) + \alpha(G)$$ 
$$= \chi(G) + \omega(G^c) \leq \chi(G) + \chi(G^c) \leq n + 1.$$
Then $\chi(G) = n - k$ and $\alpha(G) = \omega(G^c) = \chi(G^c) = k + 1$. Next we will present 
a claim and its proofs. \\

\noindent {\bf Claim 1.}  $G[V(G) - S]$ is a complete. \\
 
\noindent{\bf Proof of Claim 1.} Suppose, to the contrary, that $G[V(G) - S]$ is not complete. Then there exist vertices 
$x$, $y \in V(G) - S$ such that $xy \not \in E(G)$. Notice that $S$ is independent in $G$. We can have a proper coloring for $G$ in the following way. 
Use $n - |S| - 1 = n - k - 2$ different colors to color the vertices 
in $V(G) - S$ and another color to color all the vertices in $S$.  Therefore
$n - k = \chi(G) \leq |V(G) - S| - 1 + 1 = n - k - 1$, a contradiction. \hfill{$\diamond$} \\

Set $T_i := \overrightarrow{C}[u_i^{++}, u_{i + 1}]$, where $1 \leq i \leq k$ and the index $k + 1$ is regarded as $1$.
Obviously, $|T_i| \geq 1$ for each $i$ with $1 \leq i \leq k$.   
Set $T := \{\, i : |T_i| \geq 2 \,\}$. Next we, according to the different sizes of $|T|$, divide the remainder of the proofs into three cases. \\

{\bf Case 0} $\,\,|T| = 0.$ \\

Since $|T| = 0$, we have $C = u_1 u_1^+ u_2 u_2^+ ... u_k u_k^+ u_1$. 
We first consider the case of $|V(G) - V(C)| \geq  2$. 
Since $\alpha(G) = k + 1$, we have, by Lemma $2$, that $G[V(G) - V(C)]$ is complete.
Let $z$ be a vertex in $V(G) - V(C) - \{\, x_0\,\}$. Then $x_0 z\in E(G)$. Since $z \in V(G) - S$ and $u_2 \in V(G) - S$, we, by Claim $1$,
have that $zu_2 \in E(G)$. Thus  $G$ has a cycle $x_0 z\overrightarrow{C}[u_2, u_1]P_1x_0$ 
which is longer than $C$, a contradiction. \\

Next we consider the case $V(G) - V(C) = \{\, x_0\,\}$.
Since $S$ is independent and  $d(w) \geq \delta \geq k$ for each vertex $w$ in $S$, we must have that
for any vertex $x \in S$ and any vertex $y \in V(G) - S$, $xy \in E(G)$. Thus $G$ is $K_k \vee K^c_{k + 1}$. 
Namely, $G$ is $K_k \vee (K_k^c \cup K_{n - 2k})$ with $n = 2 k + 1$\\

{\bf Case 1} $\,\,|T| = 1.$ \\

Without loss of generality, we assume that $|T_1| \geq 2$,  $|T_r| = 1$ for each $r$ with $2 \leq r \leq k$. 
We first consider the case of $|V(G) - V(C)| \geq  2$. 
Since $\alpha(G) = k + 1$, we have, by Lemma $2$, that $G[V(G) - V(C)]$ is complete.
Let $z$ be a vertex in $V(G) - V(C) - \{\, x_0\,\}$. Then $x_0 z\in E(G)$. Since $z \in V(G) - S$ and $u_3 \in V(G) - S$, we, by Lemma $1$,
have that $zu_3 \in E(G)$. Notice that $u_3$ is regarded as $u_1$ when $k = 2$. Thus  $G$ has a cycle $x_0 z\overrightarrow{C}[u_3, u_2]P_2x_0$ 
which is longer than $C$, a contradiction. \\

Next we consider the case $V(G) - V(C) = \{\, x_0\,\}$. Obviously, we now have that $n \geq 2 k + 2$. 
Let $T_1 = y_1y_2 ... y_ru_2$, where $r \geq 1$. We first notice that
$y_rx_0 \not \in E(G)$ and $y_r u_s^+ \not \in E(G)$, where $2 \leq s \leq k$ otherwise $G$ would have cycles which are longer than $C$. We further notice that $u_1^+y_r \in E(G)$ otherwise $\{x_0, u_1^+, u_2^+, ..., u_k^+, y_r \}$ would be an independent set of size $k + 2$. 
We claim that $x_0y_{r - 1} \not \in E(G)$ otherwise $G$ would have a cycle $x_0\overleftarrow{C}[y_{r - 1}, u_1^+]\overrightarrow{C}[y_r, u_1]P_1x_0$
which is longer than $C$. We further claim that $u_l^+ y_{r - 1} \not \in E(G)$ for each $l$ with $2 \leq l \leq k$ otherwise
$G$ would have a cycle $x_0P_l\overleftarrow{C}[u_l, y_r]\overrightarrow{C}[u_1^+, y_{r - 1}]\overrightarrow{C}[u_l^+, u_1]P_1x_0$ 
which is longer than $C$. Since 
$\{x_0, u_1^+, u_2^+, ..., u_k^+, y_{r - 1}\}$ is not independent, we must have that $u_1^+ y_{r - 1} \in E(G)$.
Repeating this process, we can prove that $y_j u_1^+ \in E(G)$ for each $j$ with $1 \leq j \leq r$,
$x_0y_j \not \in E$ for each $j$ with $1 \leq j \leq r$, 
and  $y_j u_l^+ \not \in E(G)$ for each $j$ and $l$ with $1 \leq j \leq r$ and $2 \leq l \leq k$. \\

Notice that $d(w) \geq \delta \geq k$ for each vertex $w \in S - \{\, u_1^+ \,\} = \{x_0, u_2^+, ..., $ $u_k^+ \}$. We must have that
$wu_s \in E(G)$ for each vertex $w \in S - \{\, u_1^+ \,\}$ and each $s$ with $1 \leq s \leq k$. Next we will prove that
$u_1^+u_t \in E(G)$ for each $t$ with $1 \leq t \leq k$. Obviously, $u_1^+u_1 \in E(G)$. Without loss of generality, we assume that 
$u_1^+ u_2 \not \in E(G)$. In this case, we can have a proper coloring for $G$ in the following way. 
Firstly, use $n - |S| = n - k - 1$ different colors to color the vertices 
in $V(G) - S$$= \{\, y_1, y_2, ..., y_r, u_1, u_2, ..., u_k\,\}$; secondly, use the color assigned to the vertex $u_2$ in the first step of coloring to color $u_1^+$; finally, use the
color assigned to $y_1$ in the first step of coloring to color the vertices $x_0$, $u_2^+$, $u_3^+$, ... , $u_k^+$.  Thus
$n - k = \chi(G) \leq n - k - 1$, a contradiction.\\

Now $G[V(G) - S]$ is complete, $u_1^+w \in E(G)$ for each $w \in V(G) - S$, $xy \in E(G)$ 
for each vertex $x \in \{x_0, u_2^+, ..., u_k^+ \}$ and each $y \in  \{u_1, u_2, ..., u_k \}$,
and $xy \not \in E(G)$ for each vertex $x \in \{x_0, u_2^+, ..., u_k^+ \}$ and each $y \in  \{y_1, y_2, ..., y_r \}$,
we have that $G$ is $K_k \vee (K_k^c \cup K_{n - 2k})$ with $n \geq 2 k + 2$. \\

{\bf Case 2} $\,\,|T| \geq 2.$ \\

Notice that $T_i \subseteq V(G) - S$ for each $i$ with $1 \leq i \leq k$. We, by Claim $1$, have that 
$G[T_1 \cup T_2 \cup \cdots \cup T_k]$ is complete. Since $|T| \geq 2$, 
there exist two different indexes $i$ and $j$ such that $u_i^- \in T_i$, $u_j^- \in T_j$ and 
therefore $u_i^- u_j^- \in E(G)$, where $1 \leq i, \, j \leq k$. 
Then we can easily find a cycle in $G$ which is longer than $C$, a contradiction.  \\

So the proof of Theorem $1$ is completed. \hfill{$\diamond$} \\

\end{document}